\documentclass[final]{amsart}
\usepackage{amscd}
\usepackage{amsxtra}
\usepackage{amssymb, amsmath}
\usepackage{amsthm}

\numberwithin{equation}{section}

\theoremstyle{plain}
\newtheorem{theorem}{Theorem}[subsection]
\newtheorem{proposition}[theorem]{Proposition}
\newtheorem{lemma}[theorem]{Lemma}
\newtheorem{corollary}[theorem]{Corollary}
\newtheorem{definition}[theorem]{Definition}
\newtheorem{remark}[theorem]{Remark}

\newtheorem{condition}[theorem]{Condition}
\newtheorem{maintheorem}[theorem]{Main Theorem}
\begin{document}

\newcommand{\U}{\ensuremath{\mathfrak{u}_{1}}}

\newcommand{\ad}{\ensuremath { ad(\mathfrak{u}_{1})}}
\newcommand{\Oo}{\ensuremath{\varOmega^{0}(ad(\mathfrak{u}_{1}))}}
\newcommand{\Ok}{\ensuremath {\varOmega^{1}(ad(\mathfrak{u}_{1}))}}

\newcommand{\s}{\ensuremath{\mathcal{S}}}
\newcommand{\csa }{\ensuremath {\mathcal{S_{\alpha}}}}
\newcommand{\pcsa }{\ensuremath {\mathcal{S^{+}_{\alpha}}}}
\newcommand{\cs}{\ensuremath{\mathcal{S^{c}}}}
\newcommand{\vsa}{\ensuremath{\varGamma(\mathcal{S}^{+}_{\alpha})}}
\newcommand{\sad}{\ensuremath{\varGamma^{\mathcal{D}}(\mathcal{S^{+}_{\alpha}})}}
\newcommand{\san}{\ensuremath{\varGamma^{\mathcal{N}}(\mathcal{S^{+}_{\alpha}})}}
\newcommand{\la }{\ensuremath{{\mathcal{L}}_{\alpha}}}
\newcommand{\sla }{\ensuremath{{\mathcal{L}}^{1/2}_{\alpha}}}

\newcommand{\ca}{\ensuremath{\mathcal{C}_{\alpha}}}
\newcommand{\cca}{\ensuremath{\mathcal{C}^{\mathfrak{C}}_{\alpha}}}
\newcommand{\cad}{\ensuremath{\mathcal{C}^{\mathcal{D}}_{\alpha}}}
\newcommand{\can}{\ensuremath{\mathcal{C}^{\mathcal{N}}_{\alpha}}}
\newcommand{\Aa}{\ensuremath {\mathcal{A}_{\alpha}}}
\newcommand{\Aad}{\ensuremath {\mathcal{A}^{\mathcal{D}}_{\alpha}}}
\newcommand{\Aan}{\ensuremath {\mathcal{A}^{\mathcal{N}}_{\alpha}}}
\newcommand{\Q}{\ensuremath{\mathcal{A}_{\alpha}\times_\mathcal{G_{\alpha}}
\varGamma (S^{+}_{\alpha})}}
\newcommand{\B}{\ensuremath{\mathcal{B}_{\alpha}}}

\newcommand{\G}{\ensuremath{\mathcal{G}_{\alpha}}}
\newcommand{\wG}{\ensuremath{\widehat{\mathcal{G}}_{\alpha}}}

\newcommand{\h}{\ensuremath{H^{\mathcal{S}\mathcal{W}}_{(A,\phi)}}}
\newcommand{\hh}{\ensuremath{\widehat{H}^{\mathcal{S}\mathcal{W}}_{(A,\phi)}}}

\newcommand{\lda}{\lambda}

\newcommand{\spinc}{\ensuremath{Spin^{c}_{4}}}
\newcommand{\spin }{\ensuremath{Spin_{4}}}

\newcommand{\Z }{\ensuremath{\mathbb {Z}}}
\newcommand{\R }{\ensuremath{\mathbb {R}}}
\newcommand{\C }{\ensuremath {\mathbb {C}}}
\newcommand{\N }{\ensuremath {\mathbb {N}}}

\newcommand{\iso }{\ensuremath {\thickapprox }}

\newcommand{\sw}{\ensuremath {\mathcal{S}\mathcal{W}_{\alpha}}}
\newcommand{\y}{\ensuremath {\mathcal{Y}\mathcal{M}}}
\newcommand{\yp}{\ensuremath {\mathcal{Y}\mathcal{M}^{+}}}
\newcommand{\yn}{\ensuremath {\mathcal{Y}\mathcal{M}^{-}}}
\newcommand{\Cf }{\ensuremath {\mathcal{C}_{\alpha}}}
\newcommand{\w }{\ensuremath {\omega}}
\newcommand{\cx}{\ensuremath {C_{X}}}

\title{Boundary Value Problems for the $2^{nd}$-order Seiberg-Witten Equations}
\author{ Dr. Celso M Doria }
\email{cmdoria@mtm.ufsc.br}
\keywords{connections,gauge fields,4-manifolds \\ MSC 58J05 , 58E50}
\maketitle

\begin{abstract}
It is shown that the non-homogeneous Dirichlet and Neuman  problems for the $2^{nd}$-order 
Seiberg-Witten equation admit a regular solution once the $\mathcal{H}$-condition ~\ref{C:01} is 
satisfied. The approach  consist in applying the elliptic techniques to
the  variational setting of the Seiberg-Witten equation. 
\end{abstract}

\section{\bf{Introduction}}

Let X be a compact smooth 4-manifold with non-empty boundary. In our context, the Seiberg-Witten 
equations are the $2^{nd}$-order  Euler-Lagrange equation of the functional 
defined in ~\ref{D:01}. When the boundary is empty, their variational aspects  were first studied in
~\cite{JPW96} and the topological ones in  ~\cite{Do01}. Thus, the main aim is
to obtain  the existence of a solution to the non-homogeneous equations whenever  $\partial X\ne \emptyset$. 
The non-emptyness of the boundary inflicts boundary conditions on the problem.
 Classicaly, this sort of  problem  is classified according with its boundary conditions 
 in   \emph{Dirichlet Problem} ($\mathcal{D}$) or  \emph{Neumann Problem} ($\mathcal{N}$).   

\subsection{Spin$^{c}$ Structure}
The space of $Spin^{c}$ structures on X is identified with

$$Spin^{c}(X)=\{\alpha + \beta \in H^{2}(X,\Z)\oplus H^{1}(X,\Z_{2}) \mid w_{2}(X)=
\alpha (mod\ \ 2)\}.$$ 

 \noindent For each $\alpha \in Spin^{c}(X)$, there is a representation 
$\rho_{\alpha}:SO_{4}\rightarrow \C l_{4}$, induced by a \emph{$Spin^{c}$} representation,
 and  a pair of vector bundles $(\csa^{+},\la)$ over X (see ~\cite{LM89}). Let
 $P_{SO_{4}}$ be the frame bundle of $X$, so 
\begin{itemize}
\item $\csa = P_{SO_{4}}\times_{\rho_{\alpha}} V = \csa^{+}\oplus
  \csa^{-}$.\\
 The bundle  $\csa^{+}$ is the positive complex spinors bundle (fibers are 
$Spin^{c}_{4}-modules$ isomorphic to $\C^{2}$) 
\item $\la=P_{SO_{4}}\times_{det(\alpha)} \C$.\\
 It is called the \emph{determinant line bundle} associated to the
 $Spin^{c}$-struture $\alpha$. ($c_{1}(\la)=\alpha$)
\end{itemize}

Thus, for  each $\alpha \in Spin^{c}(X)$ we associate a pair of bundles

$$\alpha \in Spin^{c}(X) \quad \rightsquigarrow \qquad(\la,\csa^{+}).$$

From now on, we  considered  on $X$ a Riemannian metric $g$  and on $\csa$  an hermitian 
structure $h$.

Let  $P_{\alpha}$ be the  $U_{1}$-principal bundle over X obtained as the frame bundle 
of $\mathcal{L}_{\alpha}$ ($c_{1}(P_{\alpha})=\alpha$). Also, we consider the adjoint bundles 

$$Ad(U_{1})=P_{U_{1}}\times_{Ad}U_{1}\quad ad(\U)=P_{U_{1}}\times_{ad}\U,$$

\noindent where $Ad(U_{1})$ is  a fiber bundle with fiber $U_{1}$, and $ad(\mathfrak{u}_{1})$ is 
a vector bundle with fiber isomorphic to the Lie Algebra $\mathfrak{u}_{1}$.

\subsection{The Main Theorem}

Let $\Aa$ be (formally) the  space of connections (covariant derivative) on
$\la$, $\vsa$ is the space of
sections of $\pcsa$ and $\G=\Gamma(Ad(U_{1}))$ is  the gauge group acting on $\Aa\times\vsa$ as follows:

\begin{equation}\label{E:gauge01}
g.(A,\phi)=(A+g^{-1}dg,g^{-1}\phi).
\end{equation}

\noindent  $\Aa$ is an afim space which vector space structure, after fixing
an origin, is isomorphic to 
the space $\Ok$ of $\ad$-valued 1-forms. Once a connection $\triangledown^{0}\in\Aa$ is fixed, 
a bijection $\Aa\leftrightarrow\Ok$ is explicited by
$\triangledown^{A}\leftrightarrow A$, where $\triangledown^{A}=\triangledown^{0}+A$.  
$\G=Map(X,U_{1})$, since $Ad(U_{1})\simeq X\times U_{1}$. The curvature of a 1-connection form 
$A\in\Ok$ is the 2-form $F_{A}=dA\in \varOmega^{2}(ad(\mathfrak{u}_{1}))$.

\begin{definition}\label{def:01}
\begin{enumerate}
\item the configuration space of the $\mathcal{D}$-problem is 
\begin{equation}\label{sp:01}
\cad = \{(A,\phi)\in\Aa\times \vsa\mid (A,\phi)\mid_{Y}\overset{\text{gauge}}{\sim}(A_{0},\phi_{0})\},
\end{equation}

\vspace{05pt}

\item the configuration space of the $\mathcal{N}$-problem is 
\begin{equation}\label{sp:02}
\can = \Aa\times\vsa
\end{equation}
\end{enumerate}
\end{definition}

Although each  boundary problem requires its own configuration space, the superscripts 
$\mathcal{D}$ and $\mathcal{N}$ will be used whenever the distintion is necessary, since most 
arguments works for both sort of problems.

\noindent The Gauge Group $\G$ action on each of the configuration space is given by ~\ref{E:gauge01}. 

The Dirichlet ($\mathcal{D}$) and Neumann ($\mathcal{N}$) boundary value problems associated to the 
$\sw$-equations  are the following: Let's consider $(\Theta,\sigma)\in\Ok\oplus\vsa$ and 
$(A_{0},\phi_{0})$ defined on the manifold $\partial X$ 
($A_{0}$ is a connection on $\la\mid_{\partial X}$, $\phi_{0}$ is a section of 
$\varGamma(\mathcal{S}^{+}_{\alpha}\mid_{\partial X})$). 
In this way, find $(A,\phi)\in\cad$ satisfying $\mathcal{D}$ and
$(A,\phi)\in\can$ satisfying $\mathcal{N}$, where

\begin{equation}\label{E:11}
\begin{aligned}
\mathcal{D} = 
\begin{cases}
d^{*}F_{A} + 4\Phi^{*}(\triangledown^{A}\phi)=\Theta,\\
\Delta_{A} \phi + \frac{\left(\mid \phi \mid^{2} + k_{g}\right)}{4}\phi = \sigma, \\
(A,\phi)\mid_{\partial X}\overset{\text{gauge}}{\sim}(A_{0},\phi_{0}),
\end{cases}
\end{aligned}\quad
\begin{aligned}
\mathcal{N}=
\begin{cases}
d^{*}F_{A} + 4\Phi^{*}(\triangledown^{A}\phi)=\Theta,\\
\Delta_{A} \phi + \frac{\left(\mid \phi \mid^{2} + k_{g}\right)}{4}\phi = \sigma,\\
i^{*}(*F_{A})=0, \triangledown^{A}_{\nu}\phi =0, \\
\end{cases}
\end{aligned}
\end{equation}

\vspace{05pt}

\noindent and
\begin{enumerate}
\item the operator $\Phi^{*}:\varOmega^{1}(\csa^{+})\rightarrow \varOmega^{1}(\mathfrak{u}_{1})$
is locally given by

\begin{equation}\label{E:12}
\Phi^{*}(\triangledown^{A}\phi)= \frac{1}{2}\triangledown^{A}(\mid\phi\mid^{2})=\sum_{i}<\triangledown^{A}_{i}\phi,\phi>\eta_{i},
\end{equation}

\noindent and $\eta=\{\eta_{i}\}$ is an orthonormal frame in $\Ok$. 

\vspace{05pt}

\item $i^{*}(*F_{A})= F_{4}$, where\\
$F_{4}=(F_{14},F_{24},F_{34},0)$ is the local representation of the $4^{th}$-component (normal to 
$\partial X$) of the 2-form 
of curvature in the local chart $(x,U)$ of $X$;\\
$x(U)=\{x=(x_{1},x_{2},x_{3},x_{4})\in\R^{4};\mid\mid x\mid\mid <\epsilon, x_{4}\ge 0\},$ and \\
$x(U\cap \partial X)\subset \{x\in x(U)\mid x_{4}=0\}$. Let 
$\{e_{1},e_{2},e_{3},e_{4}\}$ be the canonical base of $\R^{4}$, so $\nu=-e_{4}$ is the normal vector 
field along $\partial X$. 

\end{enumerate}

\begin{maintheorem}\label{T:M}
If  the pair $(\Theta,\sigma)\in L^{k,2}\oplus (L^{k,2}\cap L^{\infty})$ satisfies the 
$\mathcal{H}$-condition ~\ref{C:01}, then the problems $\mathcal{D}$ and $\mathcal{N}$ 
  admit a  $C^{r}$-regular solution $(A,\phi)$, whenever $2<k$ and $r< k$.  
\end{maintheorem}

\section{\bf{Basic Set Up}}

\subsection{Sobolev Spaces}

As a vector bundle E over (X,g) is endowed with a metric and a covariant derivative 
$\triangledown$, we  define the Sobolev norm of a section 
$\phi \in \varOmega^{0}(E)$ as 

$$\mid\mid\phi\mid\mid_{L^{k,p}}=\sum_{\mid i \mid =0}^{k}
(\int_{X}\mid\triangledown^{i}\phi\mid^{p})^{\frac{1}{p}}.$$

\noindent In this way, the  $L^{k,p}$-Sobolev Spaces of sections of E is defined as

$$L^{k,p}(E)=\{\phi\in \varOmega^{0}(E) \mid \quad \mid\mid\phi\mid\mid_{L^{k,p}} < \infty\}.$$ 

\noindent In our context, in which we fixed a connection $\triangledown^{0}$ on $\la$, a metric $g$ on $X$ and an 
hermitian structure on $\csa$, the Sobolev Spaces on which the basic setting is made  are the following;
\begin{itemize}
\item $\mathcal{A}_{\alpha}= L^{1,2}(\Omega^{1}(ad(\U))))$;
\item $\varGamma(\csa^{+})$ = $L^{1,2}(\varOmega^{0}(X,\csa^{+}))$;
\item $\ca = \mathcal{A}_{\alpha}\times \varGamma(\csa^{+})$;
\item $\G= L^{2,2}(X,U_{1})=L^{2,2}(Map(X,U_{1}))$.\\ ($\G$ is  an 
$\infty$-dimensional Lie Group which Lie algebra is $\mathfrak{g}=L^{1,2}(X,\U)$). 
\end{itemize}

\vspace{05pt}

\noindent The Sobolev spaces above induce a Sobolev structure on  $\cad$ and on $\can$. From now on,
the configuration spaces will be denoted by $\ca$ by ignoring the superscripts,  unless if it needed be.

The most basic analytical results needed to achieve the main result is the \emph{Gauge Fixing Lemma} 
 (Uhlenbeck - ~\cite{KU82}) and the estimate ~\ref{E:001}, both extended by Marini,A. ~\cite{Ma92} 
to manifolds with boundary;

\begin{lemma}\label{L:00}
(Gauge Fixing Lemma) - Every connection $\hat{A}\in\Aa$ is gauge equivalent, by a gauge 
transformation $g\in\G$  named Coulomb ($\mathfrak{C}$) gauge, to a connection $A\in\Aa$ satisfying 
\begin{enumerate}
\item $d_{\tau}^{*_{f}}A_{\tau}=0$ on $\partial X$,
\item $d^{*}A=0$ on $X$.
\item In the $\mathcal{N}$-problem, the connection $A$ satisfies $A_{\nu}=0$ ($\nu\perp\partial X$).
\end{enumerate}
\end{lemma}

\begin{corollary}\label{C:00001}
Under the hypothesis of ~\ref{L:00}, there exists a constant $K>0$ such that the connection $A$,
 gauge equivalent to  $\hat{A}$ by the Coulomb gauge, satisfies the following estimates:

\begin{equation}\label{E:001}
\mid\mid A\mid\mid_{L^{1,p}}\le K.\mid\mid F_{A}\mid\mid_{L^{p}}
\end{equation}

\end{corollary}

\noindent {\bf{notation:}} $*_{f}$ is the Hodge operator in the flat metric and the index $\tau$ 
denotes tangencial components.

\subsection{Variational Formulation}

A global formulation for problems $\mathcal{D}$ and $\mathcal{N}$ is made using  
the Seiberg-Witten functional;

\begin{definition}\label{D:01}
Let $\alpha\in Spin^{c}(X)$. The Seiberg-Witten functional \\$\sw :\ca\rightarrow\R$ is defined as 

\begin{equation}\label{E:21}
\sw(A,\phi) = \int_{X}\{\frac{1}{4}\mid F_{A}\mid^{2} + \mid \triangledown^{A} \phi
\mid^{2} + \frac{1}{8}\mid \phi\mid^{4} + \frac{k_{g}}{4}\mid\phi\mid^{2}\}dv_{g}+\pi^{2}\alpha^{2}.
\end{equation}

\vspace{05pt}

\noindent where $k_{g}$= scalar curvature of (X,g).
\end{definition}

\begin{remark}\label{R:01}
 The $\G$-action on $\ca$ has the following properties;
\begin{enumerate}
\item the $\sw$-functional is $\G$-invariant.
\item the $\G$-action on $\ca$ induces on $T\ca$ a $\G$-action  as follows: \\ 
let $(\Lambda,V)\in T_{(A,\phi)}\ca$ and $g\in\G$, 

$$g.(\Lambda,V)=(\Lambda,g^{-1}V)\in T_{g.(A,\phi)}\ca.$$

\vspace{05pt}

\noindent Consequently, $d(\sw)_{g.(A,\phi)}(g.(\Lambda,V))=d(\sw)_{(A,\phi)}(\Lambda,V)$.
\end{enumerate}
\end{remark}

\noindent The tangent bundle $T\ca$ decomposes as 

$$T\ca=\Omega^{1}(\ad)\oplus\vsa.$$

\vspace{05pt}
 
\noindent In this way, the 1-form $d\sw\in\Omega^{1}(\ca)$ admits a  decomposition 
 $d\sw = d_{1}\sw + d_{2}\sw$, where

$$d_{1}(\sw)_{(A,\phi)}:\Omega^{1}(\ad)\rightarrow\R,\quad d_{1}(\sw)_{(A,\phi)}.\Lambda = 
d(\sw)_{(A,\phi)}.(\Lambda,0)$$

$$d_{2}(\sw)_{(A,\phi)}:\vsa\rightarrow\R,\quad d_{2}(\sw)_{(A,\phi)}.V=d(\sw)_{(A,\phi)}.(0,V).$$

\vspace{05pt}

\noindent By performing the computations, we get 

\begin{enumerate}
\item for every $\Lambda\in \Aa$, 

\begin{align}\label{E:22}
 d_{1}(\sw)_{(A,\phi)}.\Lambda &=\frac{1}{4}\int_{X}Re\{<F_{A},d_{A}\Lambda>  +  
4<\triangledown^{A}(\phi),\Phi(\Lambda)>\}dx, 
\end{align}

\vspace{05pt}

\noindent where 
$\Phi:\varOmega^{1}(\mathfrak{u}_{1})\rightarrow\varOmega^{1}(\csa^{+})$  
is the linear operator $\Phi(\Lambda) = \Lambda(\phi)$, which dual is defined in ~\ref{E:12}, 

\vspace{05pt}

\item for every $V \in\vsa$, 

\begin{align}\label{E:23}
 d_{2}(\sw)_{(A,\phi)}.V &= \int_{X}Re\{<\triangledown^{A}\phi ,\triangledown^{A}V> + 
<\frac{\mid\phi\mid^{2} + k_{g}}{4}\phi, V>\}dx. 
\end{align}

\end{enumerate}

\vspace{05pt}

\noindent Therefore, by taking $supp(\Lambda)\subset int(X)$ and 
$supp(V)\subset int(X)$, we restrict to the interior of $X$, and so, the gradient 
of the $\sw$-functional at $(A,\phi)\in\ca$ is

\begin{equation}\label{E:24}
grad(\sw)(A,\phi) = (d^{*}_{A}F_{A}+4\Phi^{*}(\triangledown^{A}\phi) , 
\triangle_{A}\phi + \frac{\mid\phi\mid^{2}+ k_{g}}{4}\phi)
\end{equation}

\noindent It follows from the $\G$-action on $T\ca$ that

\begin{equation}\label{E:25}
grad(\sw)(g.(A,\phi)) = \left(d^{*}_{A}F_{A}+4\Phi^{*}(\triangledown^{A}\phi) , 
g^{-1}.(\triangle_{A}\phi + \frac{\mid\phi\mid^{2}+ k_{g}}{4}\phi)\right).
\end{equation}

\vspace{10pt}

An  important analytical aspect of the $\sw$-functional is the Coercivity
Lemma proved in ~\cite{JPW96};

\begin{lemma}\label{L:01}
 {\bf{Coercivity}} - For each $(A,\phi)\in\ca$, there exists $g\in\G$ and a constant 
$K^{(A,\phi)}_{\mathfrak{C}}>0$, where $K^{(A,\phi)}_{\mathfrak{C}}$ depends on 
$(X,g)$ and $\sw(A,\phi)$,  such that  

$$\mid\mid g.(A,\phi)\mid\mid_{L^{1,2}}< K^{(A,\phi)}_{\mathfrak{C}}.$$

\end{lemma} 
\begin{proof}
lemma 2.3 in ~\cite{JPW96}. The gauge transform is the Coulomb one given in the Gauge Fixing Lemma ~\ref{L:00}.
\end{proof} 

Considering the gauge invariance of the $\sw$-theory, and the fact that the gauge group $\G$ is 
a infinite dimensional Lie Group, we can't hope to handle the problem  in the general. 
So forth, we need  to restrict the problem to  the space 

\begin{equation}\label{E:26}
\cca=\{(A,\phi)\in\ca;\mid\mid(A,\phi)\mid\mid_{L^{1,2}} < K^{(A,\phi)}_{\mathfrak{C}}\},
\end{equation}

\vspace{05pt}

\noindent The superscript $\mathcal{D}$ and $\mathcal{N}$ are being ignored for simplicity, although
each one should be taken in account according with the problem.  These choice of spaces is a 
a property of the $\G$ action on $\ca$, it is suggested by the Gauge Fixing Lemma and  the 
Coercivity Lemma; this sort of propertie is not shared by  most actions. 

\section{\bf{Existence of a Solution}}

\subsection{ Non Homogeneous Palais-Smale Condition - $\mathcal{H}$}-

In the variational formulation, the problems $\mathcal{D}$ and $\mathcal{N}$ (~\ref{E:11}) are written  as

\begin{equation} \label{E:31}
\begin{aligned}
(\mathcal{D})= \begin{cases}
grad(\sw)(A,\phi)=(\Theta,\sigma),\\
(A,\phi)\mid_{\partial X}\overset{\text{gauge}}{\sim}(A_{0},\phi_{0}).
\end{cases}
\end{aligned}
\begin{aligned}
(\mathcal{N})= \begin{cases}
grad(\sw)(A,\phi)=(\Theta,\sigma),\\
i^{*}(*F_{A})=0, \triangledown^{A}_{n}\phi =0,
\end{cases}
\end{aligned}
\end{equation}
 
\vspace{10pt}

The equations in ~\ref{E:11}  may not admit a solution for any pair 
$(\Theta,\sigma)\in \Ok\oplus\vsa$. In finite dimension, 
if we consider a function $f:X\rightarrow\R$, the analogous question would be to 
find a point $p\in X$ such that,for a fixed vector $u$, $grad(f)(p)=u$.  
This question is more subtle if  $f$ is invariant by a Lie group action on $X$. 
Therefore, we need a premiss on the pair $(\Theta,\sigma)\in\Ok\oplus\vsa$;

\begin{condition}\label{C:01} 
{\bf{($\mathcal{H}$)}} - Let $(\Theta,\sigma)\in L^{1,2}(\Ok)\oplus 
(L^{1,2}(\vsa)\cap L^{\infty}(\vsa))$ 
be a pair such that there exists a sequence $\{(A_{n},\phi_{n})\}_{n\in\Z}\subset\cca$ (~\ref{E:26}) 
with the following properties;
\begin{enumerate}
\item $\{(A_{n},\phi_{n})\}_{n\in\Z}\subset L^{1,2}(\Aa)\times (L^{1,2}(\vsa)\cup L^{\infty}(\vsa))$ and 
there exists a constant $c_{\infty}>0$ such that, for all $n\in\Z$, $\mid\mid\phi_{n}\mid\mid_{\infty}<c_{\infty}$.
\item there exists $c\in\R$ such that, for all $n\in\Z$, $\sw(A_{n},\phi_{n}) < c$,
\item the sequence  $\{d(\sw)_{(A_{n},\phi_{n})}\}_{n\in\Z}\subset \left(L^{1,2}(\Ok)\oplus 
L^{1,2}(\vsa)\right)^{*}$, of linear functionals,  converges weakly   to

$$L_{\Theta} + L_{\sigma}:T\ca\rightarrow\R,$$

\noindent where

$$L_{\Theta}(\Lambda)=\int_{X}<\Theta,\Lambda>,\quad L_{\sigma}(V)=\int_{X}<\sigma,V>.$$

\end{enumerate}
\end{condition}

\subsection{Strong Converge of $\{(A_{n},\phi_{n})\}_{n\in\Z}$ in $L^{1,2}$} -

As an immediate consequency of  (~\ref{L:01}), the sequence $\{(A_{n},\phi_{n})\}_{n\in\Z}$ 
   given by the $\mathcal{H}$-condition converges to a pair $(A,\phi)$;
\begin{enumerate}
\item weakly in $\ca$, 
\item weakly in $L^{4}(\Aa\times\vsa)$,
\item strongly in $L^{p}(\Aa\times\vsa)$, for any $p<4$.
\end{enumerate} 

\begin{remark}
Let $\{A_{n}\}_{n\in\N}\subset L^{2}$ be a converging sequence in $L^{2}$ satisfying $d^{*}A_{n}=0$,
for all $n\in\N$, and  let $A=\lim_{n\to\infty}A_{n}\in L^{2}$. So,
$d^{*}A=0$, once

$$\mid <d^{*}A,\rho >\mid \le \mid A-A_{n}\mid_{L^{2}}.\mid d\rho\mid_{L^{2}},$$

\noindent for all $\rho\in\Omega^{0}(\ad)$. 
\end{remark}

\begin{theorem}\label{T:A}
{\bf{A}} - The limit $(A,\phi)\in L^{2}(\Aa\times\vsa)$, obtained as a limit of the sequence  
$\{(A_{n},\phi_{n})\}_{n\in\Z}$, is a weak solution of ~\ref{E:11}.
\end{theorem}
\begin{proof}
The proof goes along the same lines as in the $2^{nd}$-step in the proof of the Main Theorem in 
~\cite{JPW96}.
\begin{enumerate}
\item for every $\Lambda\in \Aa$, 

\begin{align}\label{E:32}
 d_{1}(\sw)_{(A_{n},\phi_{n})}.\Lambda &=\frac{1}{4}\int_{X}Re\{<F_{A_{n}},d_{A_{n}}\Lambda>  +  
4<\triangledown^{A_{n}}(\phi_{n}),\Phi(\Lambda)>\}dx \\
 &+ \int_{\partial X}Re\{\Lambda\wedge *F_{A_{n}}\}\label{E:33}
\end{align}

\noindent where 
\begin{enumerate}
\item $\Phi:\varOmega^{1}(\mathfrak{u}_{1})\rightarrow\varOmega^{1}(\csa^{+})$  
is the linear operator $\Phi(\Lambda) = \Lambda(\phi)$; it's dual is defined in ~\ref{E:12}, 

\noindent Assuming  $\phi\in L^{\infty}$ (~\ref{L:02}), it follows that 

$$\lim_{n\to\infty} d_{1}(\sw)_{(A_{n},\phi_{n})}.\Lambda =d_{1}(\sw)_{(A,\phi)}.\Lambda.$$ 

\noindent Therefore, $d_{1}(\sw)_{(A,\phi)}.\Lambda =\int_{X}<\Theta,\Lambda>$.

\vspace{05pt}

\item $\Lambda\wedge *F_{A} = -<\Lambda, F_{4}>dx_{1}\wedge dx_{2}\wedge dx_{3}$. 

\noindent Since the equation above is true for all $\Lambda$, let $supp(\Lambda)\subset\partial X$, 
so $F_{4}=0$ ($\Rightarrow i^{*}(*F_{A})=0$).

\end{enumerate}

\vspace{10pt}

\item for every $V\in\vsa$, 

\begin{align}\label{E:34}
 d_{2}(\sw)_{(A_{n},\phi_{n})}.V &= \int_{X}Re\{<\triangledown^{A_{n}}\phi_{n} ,\triangledown^{A_{n}}V>\} + 
<\frac{\mid\phi_{n}\mid^{2} + k_{g}}{4}\phi_{n}, V>\}dx  \\
&+ \int_{\partial X}Re\{<\triangledown^{A_{n}}_{\nu}\phi_{n},V>\}\label{E:35}.
\end{align}

\noindent Analougously, it follows  that $(A,\phi)$ is a weak solution of the equation

$$d_{2}(\sw)_{(A,\phi)}.V = \int_{X}<\sigma, V>.$$

\noindent So, in the $\mathcal{N}$-problem, $\triangledown^{A}_{\nu}\phi=0$. 

\end{enumerate}

\end{proof}

In order to pursue  the strong $L^{1,2}$-convergence  for the sequence $\{(A_{n},\phi_{n})\}_{n\in\Z}$,
 next we obtain an upper bound for $\mid\mid\phi\mid\mid_{L^{\infty}}$, whenever $(A,\phi)$ is a  weak solution .

\begin{lemma}\label{L:02}
Let $(A,\phi)$ be a solution of either $\mathcal{D}$ or $\mathcal{N}$ in ~\ref{E:11}, so
\begin{enumerate}
\item If $\sigma=0$, then there exists a constant $k_{X,g}$, depending on the
  Riemannian metric on X, such that

\begin{equation}\label{E:36}
\mid\mid\phi\mid\mid_{\infty}<k_{X,g}vol(X).
\end{equation}

\item If $\sigma\ne 0$, then there exist  constantc $c_{1}=c_{1}(X,g)$ and
  $c_{2}=c_{2}(X,g)$ such that

\begin{equation}\label{E:37}
\mid\mid\phi\mid\mid_{L^{p}} < c_{1}+c_{2}\mid\mid\sigma\mid\mid^{3}_{L^{3p}}.
\end{equation}

\noindent In particular, if $\sigma\in L^{\infty}$ then $\phi\in L^{\infty}$

\end{enumerate}

\end{lemma}
\begin{proof}
\noindent Fix $r\in\R$ and suppose that there is a  ball $B_{r^{-1}}(x_{0})$, 
around the point $x_{0}\in X$,  such that 

$$\mid\phi(x)\mid >r,\quad\forall x\in B_{r^{-1}}(x_{0}).$$

\noindent Define

$$\eta = \begin{cases}
\left(1 - \frac{r}{\mid\phi\mid}\right)\phi, \quad\text{if}\quad x\in B_{r^{-1}}(x_{0}),\\
0, \quad\text{if}\quad x\in X-B_{r^{-1}}(x_{0})
\end{cases}$$

\noindent So,

\begin{equation}\label{E:38}
\mid\eta\mid\le\mid\phi\mid
\end{equation}

$$\triangledown\eta = r\frac{<\phi,\triangledown\phi>}{\mid\phi\mid^{3}}\phi + 
(1 - \frac{r}{\mid\phi\mid})\triangledown\phi$$

$$\Rightarrow \mid\triangledown\eta\mid^{2}= r^{2}\frac{<\phi,\triangledown\phi>^{2}}{\mid\phi\mid^{4}}
+2r(1 - \frac{r}{\mid\phi\mid})\frac{<\phi,\triangledown\phi>^{2}}{\mid\phi\mid^{3}}+
(1 - \frac{r}{\mid\phi\mid})^{2}\mid\triangledown\phi\mid^{2}$$

$$\Rightarrow \mid\triangledown\eta\mid^{2}< r^{2}\frac{\mid\triangledown\phi\mid^{2}}{\mid\phi\mid^{2}} + 
2r(1 - \frac{r}{\mid\phi\mid})\frac{\mid\triangledown\phi\mid^{2}}{\mid\phi\mid} + 
(1 - \frac{r}{\mid\phi\mid})^{2}\mid\triangledown\phi\mid^{2}.$$

\noindent Since $r<\mid\phi\mid$,

\begin{equation}\label{E:39}
\mid\triangledown\eta\mid^{2} < 4\mid\triangledown\phi\mid^{2}.
\end{equation}

\noindent Hence, by ~\ref{E:38} and ~\ref{E:39}, $\eta\in L^{1,2}$.\\
The directional derivative of $\sw$ at  direction $\eta $ is given by

$$d(\sw)_{(A,\phi)}(0,\eta)= \int_{X}[<\triangledown^{A}\phi,\triangledown^{A}\eta> + 
\frac{\mid\phi\mid^{2} + k_{g}}{4}\mid\phi\mid(\mid\phi\mid - r)].$$

\noindent By   ~\ref{E:23}), 

$$\int_{X}[<\triangledown^{A}\phi,\triangledown^{A}\eta> + 
\frac{\mid\phi\mid^{2} + k_{g}}{4}\mid\phi\mid(\mid\phi\mid - r)] = 
\int_{X}<\sigma,(1-\frac{r}{\mid\phi\mid})\phi>.$$

\noindent However,

$$\int_{X}<\triangledown^{A}\phi,\triangledown^{A}\eta>= 
\int_{X}[r\frac{<\phi,\triangledown^{A}\phi>^{2}}{\mid\phi\mid^{3}} + 
(1-\frac{r}{\mid\phi\mid})\mid\triangledown\phi\mid^{2}]>0.$$

\noindent  So,

$$\int_{X}\frac{\mid\phi\mid^{2} + k_{g}}{4}\mid\phi\mid(\mid\phi\mid - r) < 
\int_{X}<\sigma,(1-\frac{r}{\mid\phi\mid})\phi> < 
\int_{X}\mid\sigma\mid(\mid\phi\mid-r).$$

\noindent Hence,

$$\int_{X}(\mid\phi\mid - r)\left(\frac{\mid\phi\mid^{2} + 
k_{g}}{4}\mid\phi\mid -\mid\sigma\mid\right)<0.$$

\noindent Since $r<\mid\phi(x)\mid$, whenever $x\in B_{r^{-1}}(x_{0})$, it follows that 

\begin{equation}\label{E:310}
(\mid\phi\mid^{2} + k_{g})\mid\phi\mid < 4\mid\sigma\mid, 
\quad\text{almost everywhere in}\quad  B_{r^{-1}}(x_{0}).
\end{equation}
 
\noindent There are two cases to be analysed independently;

\begin{enumerate}
\item $\sigma=0$.\\
\noindent In this case, we get

\begin{equation}\label{E:311}
(\mid\phi\mid^{2} + k_{g})\mid\phi\mid < 0, \quad\text{almost everywhere}.
\end{equation}
 
\noindent The scalar curvature plays a central role here: if $k_{g}\ge 0$ then
$\phi=0$; otherwise, 

$$\mid\phi\mid\le max\{0,(-k_{g})^{1/2}\}.$$

\noindent Since X is compact, we let $k_{X,g}=\max_{x\in X}\{0,[-k_{g}(x)]^{1/2})$, and so,

$$\mid\mid\phi\mid\mid_{\infty}< k_{X,g}vol(X).$$

\item Let $\sigma\ne 0$.\\
The inequality ~\ref{E:310} implies that

$$\mid \phi\mid^{3} + k_{g}\mid\phi\mid -4\mid\sigma\mid < 0\quad\text{a.e}.$$

\noindent Consider the polynomial 

$$Q_{\sigma(x)}(w)= w^{3} + k_{g}w -4\mid\sigma(x)\mid.$$

\noindent A estimate for $\mid\phi\mid$ is obtained by estimating the largest real number 
$w$ satisfying $Q_{\sigma(x)}(w)<0$. $Q_{\sigma(x)}$ being monic implies that 
$\lim_{w\rightarrow\infty}Q_{\sigma(x)}(w)=+\infty$. 
So, either $Q_{\sigma(x)}>0$, whenever $w>0$, or there exist a root $\rho\in(0,\infty)$. The first 
case would imply that 

$$Q_{\sigma(x)}(\mid\phi(x)\mid) > 0,\quad a.e.,$$ 

\noindent  contradicting ~\ref{E:310}. By the same argument, there 
exists a root $\rho\in(0,\infty)$ such that $Q_{\sigma(x)}(w)$ chances its sign in a 
neighboorhood of $\rho$. Let $\rho$ be the largest root in $(0,\infty)$ with this 
propertie. By the Corollary ~\ref{AC:01}, there exist constants $c_{1}=c_{1}(X,g)$ 
and $c_{2}$ such that

$$\mid\rho\mid < c_{1} + c_{2}\mid\sigma(x)\mid^{3}.$$

\noindent Consequently,

\begin{equation}\label{E:312}
\mid\phi(x)\mid < c_{1} + c_{2}\mid\sigma(x)\mid^{3},\quad \text{a.e. in $B_{r^{-1}}(x_{0})$}
\end{equation}

\noindent and 

\begin{equation}\label{E:313}
\mid\mid\phi\mid\mid_{L^{p}}< C_{1} + C_{2}\mid\mid\sigma\mid\mid^{3}_{L^{3p}},\quad
\text{restricted to $B_{r^{-1}}(x_{0})$}
\end{equation}

\noindent where $C_{1},C_{2}$ are constants  depending on $vol(B_{r^{-1}}(x_{0}))$.

\noindent  The inequality ~\ref{E:313} can be extended over X by using a 
 $C^{\infty}$ partition of unity. Moreover, if  $\sigma\in L^{\infty}$, then

\begin{equation}\label{E:314}
\mid\mid\phi\mid\mid_{\infty} < C_{1} + C_{2}\mid\mid\sigma\mid\mid^{3}_{\infty},
\end{equation}

\noindent where $C_{1},C_{2}$ are constants  depending on $vol(X)$.

\end{enumerate}
\end{proof} 

In ~\cite{JPW96}, it was  proved a sort of concentration lemma, which is extended as follows;

\begin{lemma}\label{L:03}
Let $\{(A_{n},\phi_{n})\}_{n\in\Z}$ be the sequence given by the $\mathcal{H}$-condition  ~\ref{C:01}. So,

$$\lim_{n\to\infty}\int_{X}<\Phi^{*}(\triangledown^{A_{n}}\phi_{n}),A_{n}-A>=0.$$

\end{lemma}
\begin{proof}

By equation ~\ref{E:12}, 

$$\lim_{n\to\infty}\int_{X}<\Phi^{*}(\triangledown^{A_{n}}\phi_{n}),A_{n}-A> = 
\lim_{n\to\infty}\int_{X}<\triangledown^{A_{n}}_{i}\phi_{n},\phi_{n}>.<\eta_{i},A_{n}-A>$$

\begin{align*}
&\lim_{n\to\infty}\int_{X}<\triangledown^{A_{n}}_{i}\phi_{n},\phi_{n}>.<\eta_{i},A_{n}-A>\quad \le\\
&\lim_{n\to\infty}\int_{X}\mid<\triangledown^{A_{n}}_{i}\phi_{n},\phi_{n}>\mid^{2}.
\int_{X}\mid<\eta_{i},A_{n}-A>\mid^{2}\quad\le \\
&\lim_{n\to\infty}\left[\int_{X}\mid\triangledown^{A_{n}}_{i}\phi_{n}\mid^{2}.\mid\phi_{n}\mid^{2}\right].
\int_{X}\mid A_{n}-A\mid^{2} \quad\le \\
&\lim_{n\to\infty}c_{\infty}.\left[\int_{X}\mid\triangledown^{A_{n}}_{i}\phi_{n}\mid^{2}\right].
\mid\mid A_{n}-A\mid\mid^{2}_{L^{2}} 
\quad\le\\
&\lim_{n\to\infty}c_{\infty}.\mid\mid\phi_{n}\mid\mid^{2}_{L^{1,2}}.\mid\mid A_{n}-A\mid\mid^{2}_{L^{2}}=0.
\end{align*}
\end{proof}

\begin{theorem}
{\bf{B}} - Let $(\Theta,\sigma)$ be a pair satisfying the $\mathcal{H}$- condition  ~\ref{C:01}. So, the sequence 
$\{(A_{n},\phi_{n})\}_{n\in\Z}$, given by ~\ref{C:01}, converges strongly to $(A,\phi)\in\ca$.
\end{theorem}
\begin{proof}
From ~\ref{T:A}, $\{(A_{n},\phi_{n})\}_{n\in\Z}$ converges 
weakly in $L^{1,2}$ to $(A,\phi)\in\ca$. The prove is splitted into 2 parts;
\begin{enumerate}
\item $\lim_{n\to\infty}\mid\mid A_{n}-A\mid\mid_{L^{1,2}}=0$. \\
Let $d^{*}:\Omega^{1}(\ad)\rightarrow\Omega^{0}(\ad)$. The operator $d:ker(d^{*})\rightarrow \Omega^{2}(\ad)$ 
being elliptic implies, by the fundamentalelliptic estimate, that

$$\mid\mid A_{n}-A\mid\mid_{L^{1,2}}\le c\mid\mid d(A_{n}-A)\mid\mid_{L^{2}} + \mid\mid A_{n}-A\mid\mid_{L^{2}}.$$ 

\noindent The first term in the right-hand-side is estimate as follows;

\begin{align*}
\mid\mid dA_{n}-dA\mid\mid^{2}_{L^{2}} &= \int_{X}<d(A_{n}-A),d(A_{n}-A)>=\\
=& \int_{X}<dA_{n},d(A_{n}-A)> - \int_{X}<dA,d(A_{n}-A)>=\\
=&  \int_{X}<d^{*}F_{A_{n}},A_{n}-A> - \int_{X}<d^{*}F_{A},A_{n}-A> =\\
=& \quad d(\sw)_{(A_{n},\phi_{n})}(A_{n}-A) - 4\int_{X}<\Phi^{*}(\triangledown^{A_{n}}\phi_{n}),A_{n}-A>- \\
& \quad d(\sw)_{(A,\phi)}(A_{n}-A) - 4\int_{X}<\Phi^{*}(\triangledown^{A}\phi),A_{n}-A>+o(1)\\
=& -4\left\{\int_{X}<\Phi^{*}(\triangledown^{A_{n}}\phi_{n}),A_{n}-A> + 
\int_{X}<\Phi^{*}(\triangledown^{A}\phi),A_{n}-A> \right\} \\ 
&+ o(1),\quad \lim_{n\to\infty}o(1)=0.
\end{align*}

\noindent So, it follows from ~\ref{L:03} that 
$\lim_{n\to\infty}\mid\mid A_{n}-A\mid\mid_{L^{1,2}}=0$, and consequently,
 $A_{n}\rightarrow A$ strongly in $L^{4}$.

\vspace{10pt}

\item $\lim_{n\to\infty}\mid\mid \phi_{n}-\phi\mid\mid_{L^{1,2}}=0$.

\begin{align}\label{E:315}
\mid\mid \triangledown^{0}\phi_{n}-\triangledown^{0}\phi\mid\mid^{2}_{L^{2}} &= 
\overbrace{\int_{X}<\triangledown^{0}\phi_{n},\triangledown^{0}(\phi_{n}-\phi)>}^{(1)} - 
\overbrace{\int_{X}<\triangledown^{0}\phi,\triangledown^{0}(\phi_{n}-\phi)>}^{(2)}
\end{align}

\noindent The term $(1)$ leads to 

\begin{align*}
&\int_{X}<\triangledown^{0}\phi_{n},\triangledown^{0}(\phi_{n}-\phi)> 
= \int_{X}<(\triangledown^{A_{n}}-A_{n})\phi_{n},(\triangledown^{A_{n}}-A_{n})(\phi_{n}-\phi)> = \\
&\int_{X}<\triangledown^{A_{n}}\phi_{n},\triangledown^{A_{n}}(\phi_{n}-\phi)> - 
\int_{X}<\triangledown^{A_{n}}\phi_{n},A_{n}(\phi_{n}-\phi)> -\\
&\int_{X}<A_{n}\phi_{n},\triangledown^{A_{n}}(\phi_{n}-\phi)> +
\int_{X}<A_{n}\phi_{n},A_{n}(\phi_{n}-\phi)>=
\end{align*}

\begin{align*}
&=\overbrace{d(\sw)_{(A_{n},\phi_{n})}(\phi_{n}-\phi) - 
\int_{X}\frac{\mid\phi_{n}\mid^{2}+k_{g}}{4}<\phi_{n},\phi_{n}-\phi>}^{(11)} -\\
& \overbrace{\int_{X}<\triangledown^{A_{n}}\phi_{n},A_{n}(\phi_{n}-\phi)>}^{(12)} - 
\overbrace{\int_{X}<A_{n}\phi_{n},\triangledown^{A_{n}}(\phi_{n}-\phi)>}^{(13)} +
\end{align*}

\begin{align*}
&+ \overbrace{\int_{X}<A_{n}\phi_{n},A_{n}(\phi_{n}-\phi)>}^{(14)}.
\end{align*}

\noindent The term $(2)$ in ~\ref{E:315} leads to similar terms named
$(21)$, $(22)$, $(23)$ and $(24)$.
\noindent Let's analyse each one of the  overbraced terms obtained above:

\begin{enumerate}
\item terms {\bf{(11)}} and {\bf{(21)}}. 

\begin{align*}
& d(\sw)_{(A_{n},\phi_{n})}(\phi_{n}-\phi) - 
\int_{X}\frac{\mid\phi_{n}\mid^{2}+k_{g}}{4}<\phi_{n},\phi_{n}-\phi> + o(1)=\\
& <\sigma, \phi_{n}-\phi> - \int_{X}\frac{\mid\phi_{n}\mid^{2}+k_{g}}{4}\mid\phi_{n}-\phi\mid^{2}
- \int_{X}\frac{\mid\phi_{n}\mid^{2}+k_{g}}{4}<\phi,\phi_{n}-\phi>+\\
& o(1)\le <\sigma,\phi_{n}-\phi> - \int_{X}\frac{\mid\phi_{n}\mid^{2}+k_{g}}{4}<\phi,\phi_{n}-\phi>+ o(1)\\
\le& \mid\mid\sigma\mid\mid^{2}_{L^{2}}.\mid\mid\phi_{n}-\phi\mid\mid^{2}_{L^{2}} +
\mid\mid\frac{\mid\phi_{n}\mid^{2}+k_{g}}{4}\mid\mid^{2}_{L^{2}}.\mid\mid\phi\mid\mid_{\infty}.
\mid\mid\phi_{n}-\phi\mid\mid^{2}_{L^{2}}+ o(1),
\end{align*}

\noindent where $\lim_{n\to\infty}o(1)=0$. By the similarity among  $(11)$ and $(21)$, 
we conclude the boundeness of term $(22)$.

\vspace{05pt}

\item terms {\bf{(12)}} and {\bf{(22)}}.

\noindent $i.$ $(12)$
\begin{align*}
&\int_{X}<\triangledown^{A_{n}}\phi_{n},A_{n}(\phi_{n}-\phi)> = \\
&\int_{X}<\triangledown^{A_{n}}\phi_{n},(A_{n}-A)(\phi_{n}-\phi)> + 
\int_{X}<\triangledown^{A_{n}}\phi_{n},A(\phi_{n}-\phi)>\\
\le& \int_{X}\mid\triangledown^{A_{n}}\phi_{n}\mid^{2}.\int_{X}\mid A_{n}-A\mid^{4}.\int_{X}\mid\phi_{n}-\phi\mid^{4}
+ \\
&\int\mid\triangledown^{A_{n}}\phi_{n}\mid^{2}.\int_{X}\mid A(\phi_{n}-\phi)\mid^{2}
\end{align*}

\noindent $ii.$ $(21)$
\begin{equation*}
\int_{X}<\triangledown^{A}\phi,A(\phi_{n}-\phi)> \quad
\le \int_{X}\mid\triangledown^{A}\phi\mid^{2}.\int_{X}\mid
A(\phi_{n}-\phi)\mid^{2}
\end{equation*}

\noindent The term $\int_{X}\mid\triangledown^{A}\phi\mid^{2}$ is bounded by  ~\ref{P:405}
 and $A\in C^{0}$ by ~\ref{Th:408}.

\vspace{05pt}

\item term {\bf{$\{(13)-(23)\}$}}.

\begin{align*}
&\int_{X}<A_{n}\phi_{n},\triangledown^{A_{n}}(\phi_{n}-\phi)> - 
\int_{X}<A\phi,\triangledown^{A}(\phi_{n}-\phi)>=\\
&\int_{X}<(A_{n}-A)\phi_{n},\triangledown^{A_{n}}(\phi_{n}-\phi)> + 
\overbrace{\int_{X}<A\phi_{n},\triangledown^{A_{n}}(\phi_{n}-\phi)>}^{(i)} -\\
&\int_{X}<(A_{n}-A)\phi,\triangledown^{A}(\phi_{n}-\phi)> - 
\overbrace{\int_{X}<A_{n}\phi,\triangledown^{A}(\phi_{n}-\phi)>}^{(ii)}=\\
\end{align*}

\noindent In each of the last two lines above, the first  terms  are bounded by\\
$\mid\mid A_{n}-A\mid\mid_{L^{4}}$, while the term $\{(i)-(ii)\}$ can be written as 

\begin{align*}
& \int_{X}<(A-A_{n})\phi_{n},\triangledown^{A_{n}}(\phi_{n}-\phi)> + 
\int_{X}<A_{n}(\phi_{n}-\phi),\triangledown^{A_{n}}(\phi_{n}-\phi)> +\\
& \int_{X}<A_{n}\phi,(\overbrace{\triangledown^{A_{n}}-\triangledown^{A}}^{(A_{n}-A)})(\phi_{n}-\phi)>.
\end{align*}

\noindent So, it is also bounded by $\mid\mid A_{n}-A\mid\mid_{L^{4}}$.

\vspace{05pt}

\item term $\{(14)-(24)\}$.

\begin{align*}
&\int_{X}<A_{n}\phi_{n},A_{n}(\phi_{n}-\phi)> - \int_{X}<A\phi,A(\phi_{n}-\phi)>=\\
&= \int_{X}<A_{n}\phi_{n},(A_{n}-A)(\phi_{n}-\phi)> +
\int_{X}<(A_{n}-A)\phi_{n},A((\phi_{n}-\phi)> + \\
&\int\mid A(\phi_{n}-\phi)\mid^{2}
 \end{align*}

\end{enumerate}

\noindent Since $A\in C^{0}$,  it follows that 
$\lim_{n\to\infty}\mid\mid A(\phi_{n}-\phi)\mid\mid^{2}=0$.

\end{enumerate}
\end{proof}


\section{\bf{Regularity of the Solution $(A,\phi)$}}

Let $\beta=\{e_{i};1\le i\le 4\}$ be a orthonormal frame fixed on $TX$ with the 
following properties; for all $i,j\in\{1,2,3,4\}$

\begin{enumerate}
\item (commuting)  $[e_{i},e_{j}]=0$,
\item $\triangledown_{e_{i}}e_{j}=0$ ($\triangledown$ = Levi-Civita connection on X).
\end{enumerate}

\noindent Let $\beta^{*}=\{dx_{1},\dots,dx_{n}\}$ be the dual frame induced on $\csa^{*}$. 
From the 2$^{nd}$-property of the frame $\beta$, it follows  that 
$\triangledown_{e_{i}}dx^{j}=0$ for all $i,j\in\{1,2,3,4\}$. For the sake of simplicity, let 
 $\triangledown^{A}_{e_{i}}=\triangledown^{A}_{i}$.
Therefore, $\triangledown^{A}:\Oo\rightarrow\Ok$ is given by 

$$\triangledown^{A}\phi=\sum_{l}(\triangledown^{A}_{l}\phi)
dx_{l}\quad\Rightarrow \quad\mid \triangledown^{A}\phi\mid^{2} = 
\sum_{l}\mid\triangledown^{A}_{l}\phi\mid^{2},$$

\noindent and

$$(\triangledown^{A})^{2}=\sum_{k,l}(\triangledown^{A}_{k}\triangledown^{A}_{l}\phi)
dx_{l}\wedge dx_{k}\quad\Rightarrow\quad 
\mid(\triangledown^{A})^{2}\mid^{2}=\sum_{k,l}\mid\triangledown^{A}_{k}
\triangledown^{A}_{l}\phi\mid^{2}.$$

\noindent In this setting, the 2-form of curvature of the connection $A$ is given by

$$(F_{A})_{kl} = F_{kl} = \triangledown^{A}_{l}\triangledown^{A}_{k} - 
\triangledown^{A}_{k}\triangledown^{A}_{l}.$$

\noindent In order to compute the operator 
$\Delta_{A}=(\triangledown^{A})^{*}\triangledown^{A}:\Omega^{0}(\pcsa)\rightarrow\Omega^{0}(\pcsa)$,
 let $*:\Omega^{i}(\csa)\rightarrow\Omega^{4-i}(\csa)$ be the Hodge operator and
consider the identity

$$(\triangledown^{A})^{*} = - *\triangledown^{A}*:\Omega^{1}(\pcsa)\rightarrow\Omega^{0}(\pcsa).$$

\noindent Hence, 

$$\Delta_{A}\phi = -\sum_{k}\triangledown^{A}_{k}\triangledown^{A}_{k}\phi.$$
 
\noindent In this way,  

\begin{equation}\label{E:41}
\begin{split}
&\mid\Delta_{A}\phi\mid^{2} = \sum_{k,l}<\triangledown^{A}_{k}\triangledown^{A}_{k}\phi,
\triangledown^{A}_{l}\triangledown^{A}_{l}\phi> =\\
&= \sum_{k,l} \left[\triangledown^{A}_{k}(<\triangledown^{A}_{k}\phi,\triangledown^{A}_{l}\triangledown^{A}_{l}\phi>) -
 <\triangledown^{A}_{k}\phi,\triangledown^{A}_{k}\triangledown^{A}_{l}\triangledown^{A}_{l}\phi>\right]=\\
&=\sum_{k,l}\left[\triangledown^{A}_{k}(<\triangledown^{A}_{k}\phi,\triangledown^{A}_{l}
\triangledown^{A}_{l}\phi>) - 
<\triangledown^{A}_{k}\phi,\triangledown^{A}_{l}\triangledown^{A}_{k}\triangledown^{A}_{l}\phi>
- <\triangledown^{A}_{k}\phi,F_{lk}\triangledown^{A}_{l}\phi> \right]\\
&= \sum_{k,l}\left[\triangledown^{A}_{k}(<\triangledown^{A}_{k}\phi,
\triangledown^{A}_{l}\triangledown^{A}_{l}\phi>) - 
\triangledown^{A}_{l}(<\triangledown^{A}_{k}\phi,\triangledown^{A}_{k}\triangledown^{A}_{l}\phi>)\right]+\\
&+ \sum_{k,l}\left[<\triangledown^{A}_{l}\triangledown^{A}_{k}\phi,\triangledown^{A}_{k}\triangledown^{A}_{l}\phi> + 
<\triangledown^{A}_{k}\phi,F_{lk}\triangledown^{A}_{l}\phi>\right] =
\end{split}
\end{equation}

\begin{equation}\label{E:42}
\begin{split}
&=\sum_{k,l}\left[\triangledown^{A}_{k}(<\triangledown^{A}_{k}\phi,
\triangledown^{A}_{l}\triangledown^{A}_{l}\phi>) - 
\triangledown^{A}_{l}(<\triangledown^{A}_{k}\phi,\triangledown^{A}_{k}\triangledown^{A}_{l}\phi>)\right]+
\sum_{k,l}\mid\triangledown^{A}_{k}\triangledown^{A}_{l}\phi\mid^{2} +\\
&+ \sum_{k,l}\left[<F_{kl}\phi,\triangledown^{A}_{k}\triangledown^{A}_{l}\phi> + 
<\triangledown^{A}_{k}\phi,F_{kl}\triangledown^{A}_{l}\phi>\right]
\end{split}
\end{equation}

\noindent and so,

\begin{equation}\label{E:43}
\begin{split}
&\mid (\triangledown^{A})^{2}\phi\mid^{2} \le \quad\mid\Delta_{A}\phi\mid^{2} + 
\sum_{k,l}\left\{\mid\triangledown^{A}_{k}(<\triangledown^{A}_{k}\phi,
\triangledown^{A}_{l}\triangledown^{A}_{l}\phi>)\mid\right\} + \\
&\sum_{k,l}\left\{\mid\triangledown^{A}_{l}(<\triangledown^{A}_{k}\phi,\triangledown^{A}_{k}\triangledown^{A}_{l}\phi>)
\mid\right\} +
\sum_{k,l}\left\{\mid<F_{kl}\phi,\triangledown^{A}_{k}\phi\triangledown^{A}_{l}\phi>\mid\right\} +\\
&\sum_{k,l}\left\{\mid<\triangledown^{A}_{k}\phi,F_{kl}\triangledown^{A}_{l}\phi>\mid\right\}
\end{split}
\end{equation}

\noindent Now, by applying the inequalities

$$\left(\sum_{i}a_{i}\right)^{r}\le K_{r}.\sum_{i}\mid a_{i}\mid^{r},\quad 
\sqrt{\sum_{i=1}^{n} a_{i}}\le \sum_{i=1}^{n}\sqrt{a_{i}}$$

\noindent  to ~\ref{E:43}, we get

\begin{align*}
&\mid (\triangledown^{A})^{2}\phi\mid^{p} \le \quad K_{p}.\mid\Delta_{A}\phi\mid^{p} + 
K_{p}.\sum_{k,l}\left\{\mid\triangledown^{A}_{k}(<\triangledown^{A}_{k}\phi,
\triangledown^{A}_{l}\triangledown^{A}_{l}\phi>)\mid^{\frac{p}{2}}\right\} +\\ 
& K_{p}\sum_{k,l}\left\{\mid\triangledown^{A}_{l}(<\triangledown^{A}_{k}\phi,
\triangledown^{A}_{k}\triangledown^{A}_{l}\phi>)\mid^{\frac{p}{2}}\right\}  + 
\sum_{k,l}\left\{\mid<F_{kl}\phi,\triangledown^{A}_{k}\phi\triangledown^{A}_{l}\phi>\mid^{\frac{p}{2}}\right\} +\\
&\sum_{k,l}\left\{\mid<\triangledown^{A}_{k}\phi,F_{kl}\triangledown^{A}_{l}\phi>\mid^{\frac{p}{2}}\right\};
\end{align*}

\noindent After integrating, it follows that

\begin{equation}\label{E:44}
\begin{split}
& k_{1}.\mid\mid (\triangledown^{A})^{2}\phi\mid\mid^{p}_{L^{p}}
\quad \le \quad\mid\mid\Delta_{A}\phi\mid\mid^{p}_{L^{p}} 
 + k_{2}.\mid\mid\triangledown^{A}\phi\mid\mid^{p}_{L^{p}} + 
+ k_{3}.\mid\mid F_{A}(\phi)\mid\mid^{p}_{L^{p}}    + \\
& + k_{4}.\mid\mid F_{A}(\triangledown^{A}\phi)\mid\mid^{p}_{L^{p}} + 
k_{5}.\sum_{k,l}\int_{x}\left\{\mid\triangledown^{A}_{k}(<\triangledown^{A}_{k}\phi,
\triangledown^{A}_{l}\triangledown^{A}_{l}\phi>)\mid^{\frac{p}{2}}\right\} + \\
& + k_{6}\sum_{k,l}\int_{X}\left\{\mid\triangledown^{A}_{l}(<\triangledown^{A}_{k}\phi,
\triangledown^{A}_{k}\triangledown^{A}_{l}\phi>)
\mid^{\frac{p}{2}}\right\}
\end{split}
\end{equation}

\noindent The boundness of the right hand side of ~\ref{E:44} results from the analysis of each term;

\begin{proposition}\label{P:405}
Let $(A,\phi)\in\ca$ be a solution of equations in (~\ref{E:11}). If
$\sigma\in L^{\infty}$, then 
\begin{enumerate}
\item $\triangledown^{A}\phi\in L^{2}$,
\item $\Delta_{A}\phi\in L^{2}$.
\end{enumerate}
\end{proposition}
\begin{proof}
\begin{enumerate}
\item $\triangledown^{A}\phi\in L^{2}$

\begin{align*}
&<\Delta_{A}\phi,\phi> + \left(\frac{\mid\phi\mid^{2} +
  k_{g}}{4}\right)\mid\phi\mid^{2}  =  <\sigma,\phi>\\
&\Rightarrow\quad\mid\triangledown^{A}\phi\mid^{2}  + \left(\frac{\mid\phi\mid^{2} + k_{g}}{4}\right)\mid\phi\mid^{2} 
 = <\sigma,\phi>\quad  \le\\
&\le\quad  \frac{1}{\epsilon^{2}}\mid\sigma\mid^{2}  + \epsilon^{2}\mid\phi\mid^{2} 
\end{align*}

\noindent Therefore,

$$ \mid\triangledown^{A}\phi\mid^{2} < \frac{1}{\epsilon^{2}}\mid\sigma\mid^{2} 
+ (\epsilon^{2} - \frac{k_{g}}{4})\mid\phi\mid^{2} -  \frac{\mid\phi\mid^{4}}{4}$$

\noindent From ~\ref{L:02}, there exists a polynomial $p$, which coefficients depend 
on $(X,g)$ and $\epsilon$, such that 

\begin{equation}\label{E:45}
\mid\mid\triangledown^{A}\phi\mid\mid^{2}_{L^{2}} < p(\mid\mid\sigma\mid\mid_{\infty})
\end{equation}

\noindent So,  $\triangledown^{A}\phi\in L^{2}$.

\item $\Delta_{A}\phi$ $\in$ $L^{2}$.

$$<\Delta_{A}\phi,\Delta_{A}\phi> + \frac{\mid\phi\mid^{2}  + 
 k_{g}}{4}<\phi,\Delta_{A}\phi>  = <\sigma,\Delta_{A}\phi>$$

\noindent let $0< \epsilon <1$,

\begin{align*}
\mid\Delta_{A}\phi\mid^{2} + \frac{\mid\phi\mid^{2} + k_{g}}{4}\mid\triangledown^{A}\phi\mid^{2}
 &= <\sigma,\Delta_{A}\phi>\quad < \\ 
 &< \frac{1}{\epsilon^{2}}\mid\sigma\mid^{2} + \epsilon^{2}\mid\Delta_{A}\phi\mid^{2}
\end{align*}

\begin{equation}\label{E:46}
\begin{split}
 (1-\epsilon^{2})\mid\Delta_{A}\phi\mid^{2} + \frac{\mid\phi\mid^{2} +
k_{g}}{4}\mid\triangledown^{A}\phi\mid^{2} < \frac{1}{\epsilon^{2}}\mid\sigma\mid^{2}
\end{split}
\end{equation}

\noindent By the boundness of the term

\begin{equation}\label{E:47}
\int_{X}\mid\phi\mid^{2}.\mid\triangledown^{A}\phi\mid^{2} < 
\mid\mid\phi\mid\mid^{2}_{\infty}.\mid\mid\triangledown^{A}\phi\mid\mid^{2}_{L^{2}},
\end{equation}

\noindent  it follows the existence of a polynomial $q$,
 which coefficients depending on $\epsilon$ and $(X,g)$, such that

\begin{equation}\label{E:48}
\mid\mid\Delta_{A}\phi\mid\mid_{L^{2}}< q(\mid\mid\sigma\mid\mid_{\infty}).
\end{equation}

\end{enumerate}

\end{proof}

\begin{proposition}\label{P:406}
Let $(A,\phi)$ be solutions of the $\sw$-equations where 
$(\Theta,\sigma)\in L^{1,2}\times( L^{1,2}\cap L^{\infty})$, then $F_{A}\in L^{q}$, for all $q< \infty$.
\end{proposition}
\begin{proof}
By ~\ref{E:12}, $\Phi^{*}(\triangledown^{A}\phi)=\frac{1}{2}\triangledown^{A}(\mid\phi\mid^{2})$,
 and so,

$$d^{*}F_{A} + 4\Phi^{*}(\triangledown^{A}\phi)=\Theta \quad \Rightarrow \quad
\mid\mid d^{*}F_{A}\mid\mid^{2}_{L^{2}}\le\mid\mid\phi\mid\mid^{2}_{L^{1,2}} + 
\mid\mid\Theta\mid\mid_{L^{2}}$$

\noindent There are two cases to be analysed;

\begin{enumerate}
\item $F_{A}$ is harmonic.\\
 Since the  Laplacian defined on $\U$-forms is a elliptic operator, the fundamental 
inequality for elliptic operators claims that there exists a 
constant $C_{k}$ such that 

\begin{equation}\label{E:49}
\mid\mid F_{A}\mid\mid_{L^{k+2,2}}\le \mid\mid\Delta F_{A}\mid\mid_{L^{k,2}} + 
C_{k}\mid\mid F_{A}\mid\mid_{L^{2}}.
\end{equation} 

Consequently, $F_{A}$ being harmonic implies, for all $k\in\N$, that

 $$\mid\mid F_{A}\mid\mid_{L^{k,2}}\le C_{k}\mid\mid F_{A}\mid\mid_{L^{2}},\quad 
\Rightarrow \quad F_{A}\in C^{\infty}.$$

\item $F_{A}$ is not harmonic.\\
In this case, since $\Theta\in L^{1,2}$, $\phi\in L^{\infty}$  and 

\begin{align*}
\Delta_{A}F_{A} = d(<\phi,\triangledown^{A}\phi>) + d\Theta = <\phi,F_{A}(\phi)> + d\Theta,
\end{align*}

\noindent it follows  that $F_{A}\in L^{2,2}$.

\noindent Therefore, by the Sobolev embedding theorem $F_{A}\in L^{q}$, for all $q<\infty$.
\end{enumerate}

\end{proof} 

\begin{proposition}\label{P:407}
Let $(A,\phi)$ be solutions of the $\sw$-equations where 
$(\Theta,\sigma)\in L^{1,2}\times( L^{1,2}\cap L^{\infty})$, then 
$(\triangledown^{A})^{2}\phi\in L^{p}$, for all $1<p<2$.
\end{proposition}
\begin{proof}
In ~\ref{E:44}, we must take care of the last terms;
\begin{enumerate}
\item $F(\triangledown^{A}\phi)\in L^{p}$, for all $1<p<2$.
 By Young's inequality, 

\begin{align*}
\mid\mid F(\triangledown^{A}\phi)\mid\mid_{L^{p}} \le \mid\mid F_{A}\mid\mid_{L^{\frac{2p}{2-p}}}.
\mid\mid \triangledown^{A}\phi\mid\mid_{L^{2}}.
\end{align*}

\item There is no contribution from the divergent terms, since

\begin{align*}
& \int_{x}\left\{\mid\triangledown^{A}_{k}(<\triangledown^{A}_{k}\phi,
\triangledown^{A}_{l}\triangledown^{A}_{l}\phi>)\mid^{\frac{p}{2}}\right\} \le 
\left[vol(X)\right]^{\frac{2-p}{p}}\int_{x}\left\{\mid\triangledown^{A}_{k}(<\triangledown^{A}_{k}\phi,
\triangledown^{A}_{l}\triangledown^{A}_{l}\phi>)\mid\right\}.
\end{align*}

\noindent In the same way,

\begin{align*}
&\sum_{k,l}\int_{x}\left\{\mid\triangledown^{A}_{k}(<\triangledown^{A}_{k}\phi,
\triangledown^{A}_{l}\triangledown^{A}_{l}\phi>)\mid^{\frac{p}{2}}\right\} =0\\
&\sum_{k,l} \int_{X}\left\{\mid\triangledown^{A}_{l}(<\triangledown^{A}_{k}\phi,
\triangledown^{A}_{k}\triangledown^{A}_{l}\phi>)\mid^{\frac{p}{2}}\right\} =0.
\end{align*}

\noindent The estimates above applied to ~\ref{E:44} implies that

\begin{align*}
\mid\mid(\triangledown^{A})^{2}\phi\mid\mid_{L^{p}} &\le k_{1}\mid\mid\Delta_{A}\phi\mid\mid^{p}_{L^{p}} 
 + k_{2}\mid\mid\triangledown^{A}\phi\mid\mid^{p}_{L^{p}} + k_{3}\mid\mid\triangledown^{A}\phi\mid\mid^{p}_{L^{p}} +\\
&+ k_{4}\mid\mid F_{A}(\phi)\mid\mid^{p}_{L^{p}} + 
k_{5}\mid\mid F_{A}\mid\mid_{L^{\frac{p}{2-p}}}.\mid\mid\triangledown^{A}\phi\mid\mid^{p}_{L^{p}}
\end{align*}

\end{enumerate}

\end{proof}

\noindent Thus,  $\phi\in L^{2,p}$, for all $1<p<2$. Considering that $\sigma\in L^{1,2}$,  
 the bootstrap argument applied on ~\ref{E:11} implies  that $\phi\in L^{3,p}$, for every $k\ge 2$ and 
$1<p<2$. Hence, by Sobolev embedding theorem, $\phi\in C^{0}$.

\begin{theorem}\label{Th:408}
Let $(A,\phi)$ be a solution of the $\sw$-equations where 
$(\Theta,\sigma)\in L^{k,2}(\Ok)\oplus( L^{k,2}(\vsa)\cap L^{\infty}(\vsa))$, then 
$(A,\phi)\in L^{k+2,p}\times(L^{k+2,2}\cap L^{\infty})$, for all $1<p<2$. Moreover, if $k>2$, then 
$(A,\phi)\in C^{r}\times C^{r}$, for all $r<k$. 
\end{theorem}
\begin{proof}
\begin{enumerate}
\item If $\Theta\in L^{k,2}$, then by ~\ref{P:406} $F_{A}\in L^{k+1,2}$. Consequently, by ~\ref{C:00001},
$A\in L^{k+2,2}$.
\item The Sobolev class of $\phi$ is obtained by the bootstrap argument.
\end{enumerate}
\end{proof}

\appendix
 
\section{Estimates for solutions of  $3^{rd}$-degree equation}

Let $p,q\in\R$ and consider the equation

\begin{equation}\label{A:01}
x^{3}+px+q=0
\end{equation}

\begin{proposition}
 The solutions of (~\ref{A:01}) are given in  ~\cite{AG} by

\begin{equation}\label{A:02}
x_{1}=z_{1}+z_{2}, \quad x_{2}=z_{1}+ \lambda z_{2}\quad \text{and}\quad y_{3}=z_{1} + 
\lambda^{2}z_{2},
\end{equation}

\noindent where 

$$z_{1}=\sqrt[3]{-\frac{q}{2}+\sqrt[2]{D}}\quad z_{2}=\sqrt[3]{-\frac{q}{2}-\sqrt[2]{D}},$$

$$D=\frac{p^{3}}{27}+\frac{q^{2}}{4},$$

\noindent and $\lambda\in\C$ satisfies $\lambda^{3}=1$. 

\begin{corollary} \label{AC:01}
 Let $q,p\in\R$ such that $q<0$ and $p<0$. So, the solutions of equation (~\ref{A:01})
are estimates according with the following cases;
\begin{enumerate}
\item $D\ge 0$

\begin{equation}\label{A:03}
\mid x_{i}\mid \le \frac{8}{3} + \frac{1}{3}\mid
  q\mid + \frac{1}{12}q^{2} + \frac{1}{81}p^{3}
\end{equation}

\item $D<0$
\begin{equation}\label{A:04}
\mid x_{i}\mid\le 3 + \frac{1}{6}q^{2} + \frac{1}{81}\mid p\mid^{3}
\end{equation}

\end{enumerate}
\end{corollary}
\begin{proof}
Since

$$\mid x_{i}\mid\le \mid z_{1}\mid + \mid z_{2}\mid$$

\noindent it is enough to estimate $\mid z_{1}\mid$ and $\mid z_{2}\mid$. The basics identity 
needed are the following: suppose $x\ge 0$, so

$$\sqrt[2]{x}\le 1+\frac{1}{2}x$$

$$\sqrt[3]{x}\le 1 + \frac{1}{3}x$$

\begin{enumerate}
\item $D\ge 0$\\
In this case, $z_{1},z_{2}\in\R$ and 

$$\mid z_{1}\mid =\sqrt[3]{\mid -\frac{q}{2} + \sqrt[2]{D}\mid} \le 1 + 
\frac{1}{3}\mid -\frac{q}{2} + \sqrt[2]{D}\mid \le 
\frac{4}{3} + \frac{1}{6}\mid q\mid + \frac{1}{6}D$$

So,

$$\mid z_{1}\mid\le \frac{4}{3} + \frac{1}{6}\mid q\mid + \frac{1}{24}q^{2} + 
\frac{1}{162}p^{3}$$

The same estimate can be obtained for $\mid z_{2}\mid$. Hence,

$$\mid x_{i}\mid \le \frac{8}{3} + \frac{1}{3}\mid q\mid + \frac{1}{12}q^{2} + 
\frac{1}{81}p^{3}$$

\item $D\le 0$\\

In this case, $z_{1},z_{2}\in \C-\R$. Since $D\in\R$, we can write  
$\sqrt[2]{D}=i\sqrt[2]{\mid D\mid}$ and

$$z_{1}=\sqrt[3]{-\frac{1}{2}q + i\sqrt[2]{D}},\quad 
z_{2}=\sqrt[3]{-\frac{1}{2}q - i\sqrt[2]{D}}$$

Therefore,

$$\mid z_{i}\mid^{2}=\sqrt[3]{\frac{q^{2}}{4} + \mid D\mid} < 1 + \frac{1}{12}q^{2} + 
\frac{1}{3}\mid D\mid \le 1 +\frac{1}{6}q^{2} + \frac{1}{81}\mid p\mid^{3}$$

and

$$\mid z_{i}\mid < \frac{3}{2} + \frac{1}{12}q^{2} + \frac{1}{162}\mid p\mid^{3} $$

Hence,

$$\mid x_{i}\mid < 3 + \frac{1}{6}q^{2} + \frac{1}{81}\mid p\mid^{3} $$

\end{enumerate}
\end{proof}

\end{proposition}

\vspace{20pt}

\noindent\emph{Universidade Federal de Santa Catarina \\
    Campus Universitario , Trindade \\ Florian\'opolis - SC , Brasil\\
               CEP: 88.040-900 }

\noindent{http://www.mtm.ufsc.br}

\vspace{10pt}


\begin{thebibliography}{15}


\bibitem{AJ78}
ATIYAH, M; JONES, J.: \emph{Topological Aspects of Yang-Mills Theory, 
Comm.Math.Physics} {\bf{61}} (1979), 97-118.

\bibitem{Do01}
DORIA, C.M.: \emph{The Homotopy Type of the Seiberg-Witten Configuration Space},
 pre-print 2001. 

\bibitem{Do96}
DONALDSON, S.K.: \emph{The Seiberg-Witten Equations and 4-Manifold Topology},
Bull.Am.Math.Soc., New Ser. {\bf{33}}, n$^{o}$1 (1996), 45-70.


\bibitem{DK91}
DONALDSON, S.K.; KRONHEIMER, P.: \emph{The Geometry of  4-Manifold}, Oxford
University Press, 1991.

\bibitem{EL80}
EELLS, J.; LEMAIRE, L.: \emph{Selected Topics in Harmonic Maps}, CBMS {\bf{n$^{o}$50}},
AMS, 1980.

\bibitem{FU84}
FREED, D.; UHLENBECK, K.: \emph{Instantons and Four Manifolds}, MSRI
Publications, Vol 1, Springer-Verlag, 1984

\bibitem{GT83}
GILBARG, D.; TRUDINGER, N.S.: \emph{Elliptic Partial Differential Equations of Second Order},
 2$^{nd}$-edition, SCSM 224, Springer-Verlag, 1983.

\bibitem{AG}
GON\c CALVES, A.; \emph{Introdu\c c\~ao \`a \'Algebra} , Projeto Euclides, SBM

\bibitem{JPW96}
JOST, J.; PENG, X.; WANG, G.: \emph{Variational Aspects of the Seiberg-Witten
  Functional}, Calculus of  Variation {\bf{4}} (1996) , 205-218.

\bibitem{Ka84}
KATO, T.: \emph{Perturbation Theory for Linear Operators}, 2$^{nd}$-edition, SCSM 132, 
Springer-Verlag, 1984

\bibitem{LM89}
LAWSON, H.B.; MICHELSON, M.L.: \emph{Spin Geometry}, Princeton University Press, 1989.

\bibitem{Ma92}
MARINI, A.; \emph{Dirichlet and Neumenn Boundary Value Problems for Yang-Mills Connections},
Commuc. on Pure and Applied Math, {\bf{XLV}} (1992), 1015-1050.


\bibitem{JM96}
MORGAN, J.: \emph{The Seiberg-Witten Equations and Applications to the Topology 
of Smooth Four-Manifolds}, Math. Notes {\bf{44}}, Princeton Press.

\bibitem{Pa68}
PALAIS, R.S.: \emph{Foundations of Global Non-Linear Analysis}, Benjamin, inc, 1968.

\bibitem{KU82}
UHLENBECK, K.: \emph{Connections with $L^{p}$ bounds on Curvature}, Comm. Math. Phys. {\bf{83}}, 
1982, pp 31-42.







\end{thebibliography}
\end{document}